\documentclass[a4paper, 10pt]{article}

\usepackage{amsfonts}

 \usepackage{mathrsfs,amsfonts,amsmath}
 \usepackage{color}

 \makeatletter
 \@addtoreset{equation}{section}
 \makeatother
 \newtheorem{theorem}{Theorem}[section]
 \newtheorem{definition}{Definition}[section]
 \newtheorem{hypothesis}{Hypothesis}[section]
 \newtheorem{lemma}{Lemma}[section]
 \newtheorem{proposition}{Proposition}[section]
 \newtheorem{corollary}{Corollary}[section]
 \newtheorem{remark}{Remark}[section]
 \newtheorem{example}{Example}[section]

 \def\beqlb{\begin{eqnarray}}\def\eeqlb{\end{eqnarray}}
 \def\beqnn{\begin{eqnarray*}}\def\eeqnn{\end{eqnarray*}}

 \def\mbb{\mathbb}

 \def\proof{\noindent{\it Proof.~~}}
 \def\qed{\hfill$\Box$\medskip}

\newcommand{\bcen}{\begin{center}}
\newcommand{\ecen}{\end{center}}
\newcommand{\bgeqn}{\begin{equation}}
\newcommand{\edeqn}{\end{equation}}

\def\dz{\delta}
\def\lz{\lambda}

\def\ez{\epsilon}
\def\rar{\rightarrow}
\def\l{\left}
\def\r{\right}
\def\lc{\lceil}
\def\rc{\rceil}

\newcommand{\mP}{\mathbb{P}}
\renewcommand{\P}{\mathbb{P}}

\newcommand{\mE}{\mathbb{E}}
\newcommand{\W}{{\cal W}}
\begin{document}

\title{\LARGE\textbf{ A note on the scaling limits of contour functions of Galton-Watson trees}}

\maketitle

\bigskip

\centerline{Hui He and Nana Luan}

{\narrower{\narrower{\narrower

\begin{abstract}
  Recently, Abraham and Delmas constructed the distributions of super-critical L\'evy trees truncated at a fixed height by connecting super-critical L\'evy trees to (sub)critical L\'evy trees via a martingale transformation. A similar relationship also holds for discrete Galton-Watson trees. In this work, using the existing works on the convergence of contour functions of (sub)critical trees, we prove that the contour functions of truncated super-critical Galton-Watson trees converge weakly to the distributions constructed by Abraham and Delmas.
\end{abstract}

\bigskip

\bigskip


\bigskip


\par}\par}\par}

\section{Introduction}
In this note, we are interested in studying the scaling limits of contour functions of Galton-Watson trees. Since the (sub)critical case has been extensively studied, see e.g. \cite{[DL02]} and \cite{[K12]}, we mainly focus on the super-critical case.

In \cite{[DL02]}, it was shown that the scaling limits of contour functions of (sub)critical Galton-Watson trees are the  height processes  which encode the L\'evy trees;  see also \cite{[K12]} and references therein for some recent developments.  In the super-critical case, however, it is not so convenient to use contour functions to characterize either Galton-Watson trees or L\'evy trees which may have infinite mass.
To this end, Abraham and Delmas in \cite{ad:ctvmp}  showed that the distributions of super-critical continuous state branching processes (`CSBPs' in short) stopped at a fixed time are absolutely continuous w.r.t. (sub)critical CSBPs, via a martingale transformation. Since the L\'evy trees code the genealogy of CSBPs, they further defined the distributions of the super-critical L\'evy tree truncated at a fixed height by a similar change of probability.

In this work, we shall show that such distributions defined in \cite{ad:ctvmp} arise as the weak limits of scaled contour functions of super-critical Galton-Watson trees cut at a given level. Our main result is Theorem \ref{Main}. A key to this result is the observation shown in Lemma \ref{Lem:GWGir} which could be regarded as a discrete counterpart of martingale transformation for L\'evy trees constructed in \cite{ad:ctvmp}. Then by a collection of related results on convergence of Galton-Watson processes to CSBPs, we obtain our main result.

Let us mention some related works here. In \cite{[BPS12]}, the authors studied the local time processes of contour functions of binary Galton-Watson trees in continuous time. Duquesne and Winkel in \cite{[DW07]} and \cite{[DW12]} also constructed super-critical L\'evy trees as increasing limits of Galton-Watson trees. In their works, the trees are viewed as metric spaces and the convergence holds in the sense of the Gromov-Hausdorff distance; see Remark \ref{remDW} below for a discussion of the relation between the present work and \cite{[DW12]}.

The paper is organized as follows. In Section \ref{Sectree}, we recall some basic definitions and facts on trees and branching processes.  In Section \ref{Secmain}, we present our main result and its proof.

 We shall assume  that all random variables in the paper are defined on the same probability space $(\Omega, {\cal F}, {\mP}).$  Define ${\mbb N}=\{0,1,2,\ldots\}$, ${\mbb R}=(-\infty,\infty)$ and ${\mbb R}^+=[0,\infty).$

\section{Trees and branching processes}\label{Sectree}

\subsection{Discrete trees and Galton-Watson trees}\label{SecGW}
We present the framework developed in \cite{[Ne86]} for trees; see
also \cite{[LeG05]} for more notation and terminology. Introduce the set of labels
$${\mbb U}=\bigcup_{n=0}^{\infty}(\mathbb{N}^*)^n,
$$
where $\mbb{N}^*=\{1,2,\ldots\}$ and by convention
$(\mbb{N}^*)^0=\{\emptyset\}.$

An element of $\mbb U$ is thus a sequence $w = (w^1, \ldots , w^n)$
of elements of $\mbb{N}^*$, and we set $|w|=n$, so that $|w|$
represents the generation of $w$ or the height of $w$. If $w = (w^1,
\ldots, w^m)$ and $v = (v^1,\ldots, v^n)$ belong to $\mbb U$, we
write $wv = (w^1,\ldots,w^m, v^1,\ldots, v^n)$ for the concatenation
of $w$ and $v$. In particular $w\emptyset =\emptyset w = w$.

A (finite or infinite) rooted ordered tree $\mathbf t$ is a subset of
$\mbb U$ such that
\begin{enumerate}
\item $\emptyset\in {\mathbf t}.$

\item $(w^1, \ldots , w^n)\in {\mathbf t}\setminus\{\emptyset\}\Longrightarrow (w^1, \ldots , w^{n-1})\in {\mathbf t}.$

\item For every $w\in {\mathbf t}$, there exists a finite integer $k_{w}{\mathbf t}\geq
0$ such that if $k_{w}{\mathbf t}\geq1$, then $wj\in {\mathbf t}$ for any $1\leq
j\leq k_w\mathbf t$ ($k_w{\mathbf t}$ is the number of children of $w\in{\mathbf
t}$).
\end{enumerate}

Then $\emptyset$ is called the root of tree ${\bf t}$. Let ${\mathbf T}^{\infty}$ denote the set of all such trees $\mathbf t$. For each $u\in{\mbb U}$, define ${\bf T}_{u}=\{{\bf t}\in{\bf T}^{\infty}: u\in {\bf t}\}$. We endow ${\bf T}^{\infty}$ with the $\sigma$-algebra $\sigma({\bf T}_{u}, u\in{\mbb U})$; see \cite{[Ne86]} for details.
Given a tree $\mathbf t$, we call an element in the set ${\mathbf
t}\subset{\mbb U}$ a node of ${\mathbf t}$.  Denote the height of a tree
$\mathbf t$ by $|{\mathbf t}|:=\max\{|\nu|:\nu\in{\mathbf t}\}$. For $h=0,1,2,\ldots,$
define $r_h{\mathbf t}=\{\nu\in{\mathbf t}:|\nu|\leq h\}$, which is a finite tree.
Denote by $\#{\mathbf t}$ the number of nodes of  ${\mathbf
t}$. Let
$$
{\mathbf T}:=\{{\mathbf t}\in{\mathbf T}^{\infty}: \#{\mathbf t}<\infty \}
$$
be the set of all finite trees.

We say that $w\in{\mathbf t}$ is a
leaf of ${\mathbf t}$ if $k_w{\mathbf t}=0$ and set
$${\cal L}({\mathbf t}):=\{w\in{\mathbf t}: k_w{\mathbf t}=0\}.$$ So ${\cal L}({\mathbf t})$  denotes
the set of leaves of ${\mathbf t}$ and  $\#{\cal L}(\mathbf t)$ is the
number of leaves of ${\mathbf t}$.

\medskip

 To code the finite trees, we introduce the so-called {\em contour functions}; see \cite{[LeG05]} for details.  Suppose that the tree $\bf t\in \bf T$ is embedded in the half-plane in such a way that edges have length one. Imagine that  a particle starts at time $s=0$ from the root of the tree and then explores the tree from the left to the right, moving continuously
along the edges at unit speed, until all edges have
been explored and the particle has come back to the root. Then the total time needed to explore the tree is $\zeta({\bf t}) := 2(\#{\bf t}-1)$.
The contour function of $\bf t$ is the function $(C({\bf t}, s),0\leq s\leq\zeta({\bf t}))$ whose value at time $s\in [0,\zeta({\bf t})]$ is the
distance (on the tree) between the position of the particle at time $s$ and the root. We set $C({\bf t}, s)=0$ for $s\in[\zeta({\bf t}), 2\#{\bf t}].$

\medskip

Given a probability distribution $p=\{p_n: n=0,1,\ldots\}$ with
$p_1<1$, following \cite{[Ne86]} and \cite{[AP98]}, call a ${\mathbf T}^{\infty}$-valued random variable ${\cal G}^p
$ a Galton-Watson tree with offspring distribution
$p$ if
\begin{enumerate}
 \item[i)] the number of children
 of $\emptyset$ has distribution $p$:
$$
\P(k_{\emptyset}{\cal G}^p=n)=p_n,\quad \forall\, n\geq0;
$$
\item[ii)]
 for each $h=1,2,\ldots,$ conditionally given $r_h{\cal G}^p={\mathbf
t}$ with $|{\mathbf t}|\leq h$,  for $\nu\in\mathbf{t}$ with $|\nu|=h$,
$k_{\nu}{\cal G}^p$ are i.i.d. random variables distributed according to $p$.
 \end{enumerate}
 The second property is called branching property. From the definition we see for ${\mathbf t}\in {\mathbf T}$
$$
\P(r_{h+1}{\cal G}^p=r_{h+1}{\mathbf t}\bigm|
r_h\mathcal{G}^p=r_h\mathbf{t})=\prod_{\nu\in r_h\mathbf{t}\setminus r_{h-1}\mathbf{t}}p_{{k_{\nu}{\mathbf
t}}},
$$
where the product is over all nodes $\nu$ of $\mathbf t$ of height $h$. We
have then
\bgeqn\label{ForGuT}
\P({\cal G}^p={\mathbf
t})=\prod_{\nu\in{\mathbf t}}p_{k_{\nu}\mathbf{t}},\quad \mathbf{t}\in \mathbf{T},
\edeqn
where the product is over all nodes $\nu$ of $\mathbf t$. Define $m(p)=\sum_{k\geq0}kp_k$. We say ${\cal G}^p$ is sub-critical (resp. critical, super-critical) if $m(p)<1$(resp. $m(p)=1, m(p)>1$).

\medskip

Given a tree $\bf t\in T^{\infty}$, for $a\in {\mbb N}$, define $r_a{\bf t}=\{\nu\in{\bf
t}:|\nu|\leq a\}$. Then $r_a{\bf t}$ is a finite tree whose
contour function is denoted by $\{C_a({\bf t}, s):s\geq0\}$.  For $k\geq0$, we denote by $Y_k({\bf t})$ the number of individuals in generation $k$: $$Y_k({\bf t})=\#\{\nu\in{\bf t}: |\nu|=k\}, \quad k\geq 0.$$

Given a probability measure $p=\{p_k: k\geq0\}$ with
$\sum_{k\geq0}kp_k >1.$  Let ${\cal G}^p$ be a super-critical Galton-Watson tree with offspring distribution $p$.
 Then ${\P}(\#{\cal G}^p<\infty)=f(p)$, where $f(p)$ is the minimal solution of the following equation of $s$:
$$
g^p(s):=\sum_{k\geq0}s^kp_k=s, \quad 0\leq s\leq1.
$$
Let $q=\{q_k: k\geq0\}$ be another probability distribution such that
$$q_k=f(p)^{k-1}p_k, \text{ for } k\geq1,\text{ and } q_0=1-\sum_{k\geq1}q_k.$$
Then $\sum_{k\geq0}kq_k<1$.
Let ${\cal G}^q$ be a subcritical GW tree
with offspring distribution $q$. Note that $\l(Y_k({\cal G}^q), k\geq0\r)$ is a Galton-Watson process starting from a single ancestor with offspring distribution $q$. We first present a simple lemma.

\begin{lemma}\label{Lem:GWGir} Let $F$ be any nonnegative measurable function on $\bf T$. Then for ${\bf t}\in {\bf T}$,
\beqlb\label{GWGirb}
{\mP}\l[{\cal G}^p={\bf t}\r]=f(p){\mP}\left[{\cal
G}^q={\bf t}\right]
\eeqlb
and for any $a\in \mbb N$,
\beqlb\label{GWGira}
{\mE}\l[F(r_a{\cal G}^p)\r]={\mE}\left[f(p)^{1-Y_{a}({\cal G}^q)}F(r_a{\cal
G}^q)\right].
\eeqlb
\end{lemma}
\proof (\ref{GWGirb}) is just  (4.8) in \cite{adh:pgwttvmp}. The proof of (\ref{GWGira}) is straightforward. Set $\text{gen}(a, {\bf t})=\{\nu\in{\bf t}: |\nu|=a\}$.
By (\ref{ForGuT}), for $\bf t\in T$, we have
$$\P(r_a{\cal G}^p=r_a{\mathbf t})=\prod_{\nu\in r_a{\mathbf t}\setminus{\cal L}(r_a{\mathbf t})}p_{k_{\nu}{\mathbf t}}\cdot \prod_{\nu\in{\cal L}(r_a{\mathbf t})\setminus \text{gen}(a, {\bf t})}p_0 $$
and then
$$
\P(r_a{\cal G}^p=r_a{\mathbf t})=f(p)^{-\left(\sum_{\nu\in r_a{\mathbf t}\setminus{\cal
L}(r_a{\mathbf t})}(k_{\nu}{\mathbf
t}-1)\right)}\left(\frac{p_0}{q_0}\right)^{\#{\cal L}(r_a{\mathbf
t})-Y_{a}({\bf t})}\P(r_a{\cal G}^q=r_a{\mathbf t}).
$$
We also have
\begin{equation}\label{eq:pzero}
q_0=1-\sum_{k=1}^{\infty}f(p)^{k-1}p_k=1+p_0/f(p)-g^p(f(p))/f(p)=p_0/f(p).
\end{equation}
Then  (\ref{GWGira})  follows from the fact that given a tree
$\mathbf{t}\in \mathbf{T}$,
$$\#\mathcal{L}({\mathbf t})=1+\sum_{\nu\in {\mathbf t}\setminus{\cal
L}({\mathbf t})}(k_{\nu}{\mathbf t}-1).$$
 \qed

\begin{remark}\label{Rem: dis}It is easy to see that  $(f(p)^{-Y_{n}({\cal G}^q)}, n\geq0)$ is a martingale with respect to ${\cal F}_n=\sigma(r_n {\cal G}^q).$ In fact, by the branching property, we have for all $0\leq m\leq n$,
\beqlb\label{mart}
{\mbb E}\l[f(p)^{-Y_{n}({\cal G}^q)}\bigg{|}r_m{\cal G}^q\r]=\l[{\mbb E}\l[f(p)^{-Y_{n-m}({\cal G}^q)}\r]\r]^{Y_{m}({\cal G}^q)}=f(p)^{-Y_{m}({\cal G}^q)}.
\eeqlb
\end{remark}
Since   contour functions code finite trees in $\mathbf{T}$, we immediately get the following result.
\begin{corollary}\label{cor:con}
For any nonnegative measurable function $F$ on $ C({\mbb R}^+, {\mbb R}^+)$ and $a\in\mbb N$,
\beqlb\label{cor:cona}{\mE}\l[F\l( C_{a}({\cal G}^p,\cdot)\r)\r]={\mE}\left[f(p)^{1-Y_{a}({\cal G}^q)}F\l( C_{a}({\cal G}^q,\cdot)\r)\right].
\eeqlb
\end{corollary}

Lemma \ref{Lem:GWGir} could be regarded as a discrete counterpart of the martingale transformation for L\'evy trees in Section 4 of \cite{ad:ctvmp}; see also (\ref{Gircon}) below in this paper. To see this, we need to introduce  continuous state branching processes and L\'evy trees.

\subsection{Continuous State Branching Processes}\label{Secbm}
%\label{sec:exc-meas}
Let $\alpha\in \mbb R$, $\beta\ge 0$ and $\pi$ be  a $\sigma$-finite measure
on $(0,+\infty)$ such that $\int_{(0,+\infty)}(1\wedge
r^2)\pi(dr)<+\infty$.
The   branching mechanism $\psi$ with characteristics $(\alpha,\beta,
\pi)$ is
defined by:
\begin{equation}
   \label{eq:psi}
\psi(\lambda)=\alpha\lambda+\beta\lambda^2
+\int_{(0,+\infty)}\left(e^{-\lambda
  r}-1+\lambda r1_{\{r<1\}}\right)\pi(dr).
\end{equation}
A c\`ad-l\`ag $\mbb R^+$-valued Markov process $Y^{\psi,x}=(Y_t^{\psi, x}, t\geq0)$ started at $x\geq0$ is called $\psi$-continuous state branching process ($\psi$-CSBP in short) if its transition kernels satisfy
$$
E[e^{-\lambda Y^{\psi,x}_t}]=e^{-xu_t(\lambda)},\quad t\geq0,\, \lambda>0,
$$
where $u_t(\lambda)$ is the unique nonnegative solution of
$$
\frac{\partial u_t(\lz)}{\partial t}=-\psi(u_t(\lz)),\quad u_0(\lz)=\lz.
$$

\noindent $\psi$ and $Y^{\psi,x}$ are said to be sub-critical (resp. critical, super-critical) if $\psi'(0+)\in(0,+\infty)$ (resp. $\psi'(0+)=0, \psi'(0+)\in[-\infty,0)$). We say that $\psi$ and $Y^{\psi,x}$ are (sub)critical if they are critical or sub-critical.

\bigskip

\noindent In the sequel of this paper, we will assume the following assumptions on $\psi$ are in force:
\begin{enumerate}

%\item[(H1)]$\beta>0$ or $\int_{(0,1)} r \pi(dr)=+\infty $.

\item[(H1)] The Grey condition holds:
\begin{equation}
\int^{+\infty}_1\frac{d\lambda}{\psi(\lambda)}<+\infty.
\end{equation}
The  Grey  condition  is  equivalent  to the  a.s.   finiteness  of  the
extinction time  of the corresponding  CSBP. This assumption is  used to
ensure that the  corresponding height process is continuous.

%\item[(H3)] $\int_{1}^{\infty}r\pi(dr)<\infty$.

\item[(H2)]  The   branching  mechanism   $\psi$  is  conservative:   for  all
  $\varepsilon>0$,
\[
\int_{(0, \varepsilon]} \frac{d\lambda}{|\psi(\lambda)|}=+\infty.
\]
The  conservative assumption  is  equivalent to  the  finiteness of  the
corresponding CSBP at all time.

\end{enumerate}
Let us remark that (H1) implies $\beta>0$ or $\int_{(0,1)} r \pi(dr)=+\infty $. And if $\psi$ is (sub)critical, then we must have $\alpha-\int_{(1,{+\infty})}r\pi(dr)\in[0,+\infty).$ We end this subsection by collecting some results from \cite{ad:ctvmp}.

\bigskip

\noindent Let $(X_t,t\geq0)$ denote the canonical process of ${\cal D}:=D({\mbb R}^+, \mbb R)$.  Let $P_x^{\psi}$ be the probability measure on $D({\mbb R}^+, \mbb R)$ such that $P_x^{\psi}(X_0=x)=1$ and $(X_t,t\geq0)$ is a $\psi$-CSBP under $P_x^{\psi}$.

\begin{lemma}\label{AD2.2}(Lemma 2.4 in \cite{ad:ctvmp})
Assume that $\psi$ is supercritical satisfying (H1) and (H2).
Then
\begin{enumerate}

\item[(i)] $P_x^{\psi}$-a.s. $X_{\infty}=\lim_{t\rar\infty}X_t$ exists, $X_{\infty}\in\{0,\infty\}$ and
    $$
    P_x^{\psi}(X_{\infty}=0)=e^{-\gamma x},
    $$
where $\gamma$ is the largest root of $\psi(\lz)=0$.
\item[(ii)] For any nonnegative random variable measurable w.r.t. $\sigma(X_t, t\geq0)$, we have
    $$
    E_x^{\psi}[W|X_{\infty}=0]=E_x^{\psi_{\gamma}}[W],
    $$
where $\psi_{\gamma}(\cdot)=\psi(\cdot+\gamma).$

%\item[(iii)] The process $e^{\gamma X_t}$ under $P_x^{\psi_{\gamma}}$ is a %martingale w.r.t $\sigma(X_r, 0\leq r\leq t)$ and
%$$
%E_{x}^{\psi_{\gamma}}[e^{\gamma X_t}]=e^{\gamma x}.
%$$
\end{enumerate}
\end{lemma}
%\proof See Lemma 2.4 in \cite{ad:ctvmp} for (i) and (ii).
%See Theorem 2.2 in \cite{ad:ctvmp} for (iii).
%\qed
\subsection{Height processes}\label{Secbpheight}

To code the genealogy of the $\psi$-CSBP, Le Gall and Le Jan \cite{[LeLa98]} introduced the so-called height process, which is a functional of the L\'{e}vy process with Laplace exponent $\psi$; see also Duquesne and Le Gall \cite{[DL02]}.

Assume that $\psi$ is (sub)critical satisfying (H1). Let ${\mbb P}^{\psi}$ be a probability measure  on $\cal D$ such that  under ${\mbb P}^{\psi}$, $X=(X_t,t\geq0)$ is a L\'{e}vy process with
nonnegative jumps and with Laplace exponent $\psi$:
 $$
 {\mbb E}^{\psi}\Big[e^{-\lz X_t}\Big]=e^{t\psi(\lz)},\quad t\geq0,\,\lz\geq0.
 $$
%By (H1), $X$ is of infinite variation ${\mbb P}^{\psi}$-a.s.

 The so-called continuous-time \textit{height process} denoted by $H$ with sample path in $C(\mbb R^+,\mbb R^+)$ is defined for every $t\geq0$ by:
 $$
 H_t=\liminf_{\ez\rightarrow0}\frac{1}{\ez}\int_0^t1_{\{X_s<I_t^s+\ez\}}ds,
 $$
where the limit exists in ${\mbb P}^{\psi}$-probability and $I_t^s=\inf_{s\leq r\leq t}X_r$; see \cite{[DL02]}.
\noindent Under
$\mP^{\psi}$, for $a\geq0$, the local time of the height process at
level $a$ is the continuous increasing process $(L_s^a, s\geq0)$
which can be characterized via the approximation
\beqlb\label{local}
\lim_{\ez\rightarrow0}\sup_{a\geq\ez}\mE^{\psi}\left[\sup_{s\leq
t}\left|\ez^{-1}\int_0^s 1_{\{a-\ez<H_r\leq
a\}}dr-L_s^a\right|\right]=0.
\eeqlb
Furthermore, for any nonnegative measurable function $g$ on $\mbb R^+$,
\beqlb\label{localtime}
\int_0^s g(H_r)dr=\int_{\mbb R^+}g(a)L_s^ada,\quad s\geq0.
\eeqlb
%When $a>0$, we also have
%$$
%\lim_{\ez\rightarrow0}\mE^{\psi}\left[\sup_{s\geq
%t}\left[\ez^{-1}\int_0^s 1_{\{a-\ez<H_r\leq
%a\}}dr-L_s^a\right]\right]=0.
%$$

\noindent For any $x>0$, define $$T_x=\inf\{t\geq0: I_t\leq -x\},$$
 where $I_t=\inf_{0\leq r\leq t}X_r$. By Theorem 1.4.1 of \cite{[DL02]}, the process
$(L^a_{T_x},a\geq0)$ under $\mP^{\psi}$ is distributed as a $\psi$-CSBP started at $x$.

Let ${{\cal C}:=C({\mbb R}^+, {\mbb R}^+)}$ be the  space of nonnegative continuous functions on ${\mbb R}^+$ equipped with the supmum norm. Denote by $(e_t,t\geq0)$  the canonical process of ${\cal C}$. Denote by $\mP_x^{\psi}$ the law of $(H_{t\wedge T_x}, t\geq0)$ under $\mP^{\psi}$. Then $\mP_x^{\psi}$ is a probability distribution on $\cal C$. Set $Z_a=L_{T_x}^a$ under $\mP_x^{\psi}$, i.e.,
$$
\lim_{\ez\rightarrow0}\sup_{a\geq\ez}\mE_x^{\psi}\left[\sup_{s\leq
t}\left|\ez^{-1}\int_0^s 1_{\{a-\ez<e_r\leq
a\}}dr-Z_a\right|\right]=0.
$$

\subsection{Super-critical L\'{e}vy trees} \label{SecsuperCRT}
Height processes code the genealogy of (sub)critical CSBPs. However, for super-critical CSBPs, it is not so convenient to introduce the height process since the super-critical CSBPs may have infinite mass.  Abraham and Delmas in \cite{ad:ctvmp} studied the distributions of trees cut at a fixed level, where a super-critical
L\'{e}vy trees was constructed via a Girsanov transformation. We recall
their construction here.

 For any $f\in C({\mbb R}^+, {\mbb R}^+)$ and $a>0$, define
$$
\Gamma_{f,a}(x)=\int_0^x1_{\{f(t)\leq a\}}dt,\quad \Pi_{f,a}(x)=\inf\{r\geq0: \Gamma_{f,a}(r)>x\},\quad x\geq0,
$$
where we make the convention that $\inf\emptyset=+\infty.$ Then we define
$$
\pi_a(f)(x)=f(\Pi_{f,a}(x)),\quad f\in {C({\mbb R}^+, {\mbb R}^+)}, \,  x\geq0.
$$
Note that $\pi_a\circ\pi_b=\pi_a$ for $0\leq a\leq b.$ Let $\psi$ be a super-critical branching
mechanism satisfying (H2).
Denote by $q^*$  the unique (positive) root of $\psi'(q)=0$. Then the
branching mechanism $\psi_q(\cdot)=\psi(\cdot+q)-\psi(q)$ is critical for
$q=q^*$ and sub-critical for $q>q^*$. We also have $\gamma>q^*$. Because super-critical branching processes may have infinite mass,
in \cite{ad:ctvmp}  it was cut at a given level to construct the corresponding
genealogical continuum random tree.  Define
$$
M_a^{\psi_q, -q}=\exp\left\{-qZ_0+qZ_a+\psi(q)\int_0^aZ_sds\right\},\quad a\geq0.
$$
Define a filtration ${\cal H}_a=\sigma(\pi_a(e))\vee {\cal N}$, where $\cal N$ is the class of ${\mbb P}_x^{\psi_q}$ negligible sets. By (\ref{local}), we have $M^{\psi_q,-q}$ is ${\cal H}$-adapted.
\begin{theorem}(Theorem 2.2 in \cite{ad:ctvmp}) For each $q\geq q^*$, $M^{\psi_q,-q}$ is an ${\cal H}$-martingale under ${\mbb P}_x^{\psi_q}.$
\end{theorem}
\proof See Theorem 2.2 and arguments in Section 4 in \cite{ad:ctvmp}. \qed

Define the distribution $\mP_x^{\psi,a}$
%(resp.$\mN^{\psi,a}$)
of the $\psi$-CRT cut at level $a$ with initial mass $x$, as
the distribution of $\pi_a(e)$ under $M_a^{\psi_q,-q}d\mP_x^{\psi_q}$:
%(resp. $e^{qZ_a+\psi(q)\int_0^aZ_rdr}d\mN^{\psi_q}$)
for any
non-negative measurable function $F$ on $C({\mbb R}^+, {\mbb R}^+)$,
 \beqlb\label{supercriticalP}
 \mE_x^{\psi,a}[F(e)]&=&\mE_x^{\psi_q}\Big[M_a^{\psi_q,-q}F(\pi_a(e))\Big],
 %\\
 %\label{forsupercriticalN}
 %\mN^{\psi,a}[F(e)]&=&\mN^{\psi_q}
 %\Big[e^{qZ_a+\psi(q)\int_0^aZ_rdr}F(\pi_a(e))\Big],
 \eeqlb
which do not depend on the choice of $q\geq q^*$; see Lemma 4.1 of
\cite{ad:ctvmp}. Taking $q=\gamma$ in (\ref{supercriticalP}), we see
\beqlb\label{Gircon} \mE_x^{\psi,a}[F(e)]=\mE_x^{\psi_{\gamma}}\Big[e^{-\gamma x+\gamma Z_a}F(\pi_a(e))\Big]\eeqlb
and $(e^{-\gamma x+\gamma Z_a}, a\geq0)$ under ${\mP}_x^{\psi_{\gamma}}$ is an ${\cal H}$-martingale with mean 1.
\begin{remark}\label{lem: defsuper} $\mP_x^{\psi,a}$ gives the law of super-critical L\'evy trees truncated at height $a$. Then the law of the whole tree could be defined as a projective limit. To be more precise,
 let $\cal W$ be the set of $C({\mbb R}^+, {\mbb R}^+)$-valued functions endowed with the $\sigma$-field generated by the coordinate maps. Let $(w^a,a\geq0)$ be the canonical process on $\cal W$. Proposition 4.2 in \cite{ad:ctvmp}  proved that there exists a probability measure $\bar{\mP}_x^{\psi}$
%(resp.an excursion measure $\bar{\mN}^{\psi}$)
on $\W$ such that for every
$a\geq0$, the distribution of $w^a$ under $\bar{\mP}_x^{\psi}$
%(resp. $\bar{\mN}^{\psi}$)
 is $\mP_x^{\psi,a}$
 %(resp. $\mN^{\psi,a}$)
  and for $0\leq a\leq b$
$$
\pi_a(w^b)=\pi_a\quad \bar{\mP}_x^{\psi}-a.s.% (resp. \bar{\mN}^{\psi}-a.e.).
$$
\end{remark}

\begin{remark}
The above definitions of $ \mE_x^{\psi,a}$ and $\bar{\mP}_x^{\psi}$
%and $ \mN^{\psi,a}, \bar{\mN}^{\psi}$
are also valid for (sub)critical branching mechanisms.
\end{remark}

%  Let $0$ denote the path that is
%constantly zero.

\section{From Galton-Watson forests to  L\'evy forests}\label{Secmain}
%\section{Preliminary results on Galton-Watson trees and branching processes}
%\label{SecPrGW}
Comparing (\ref{cor:cona}) with (\ref{Gircon}), one can see that the l.h.s. are similar. The super-critical trees (discrete or continuum) truncated at height $a$ are connected to sub-critical trees, via a martingale transformation. Motivated by  Duquesne and Le Gall's work \cite{[DL02]}, which studied the scaling limit of (sub)critical trees, one may hope that the laws of suitably rescaled super-critical Galton-Watson trees truncated at height $a$ could converge to the law defined in (\ref{Gircon}). Our main result, Theorem \ref{Main}, will show it is true.

For each integer $n\geq1$ and real number $x>0$,
\begin{itemize}
\item
 Let $[x]$ denote the integer part of $x$ and let $\lceil x\rceil$ denote the minimal integer which is larger than $x$.

\item
 Let
$p^{(n)}=\{p^{(n)}_k: k=0,1,2,\ldots\}$ be a probability
measure on $\mbb N$.

\item Let
${\cal G}^{p^{(n)}}_1, {\cal G}^{p^{(n)}}_2, \ldots, {\cal G}^{p^{(n)}}_{[nx]}$ be independent Galton-Watson trees with the same offspring distribution $p^{(n)}.$

\item Define $Y_k^{p^{(n)},x}=\sum_{i=1}^{[nx]}Y_k({\cal G}_{i}^{p^{(n)}})$. Then $Y^{p^{(n)},x}=(Y_k^{p^{(n)},x},k=0,1,\ldots)$  is a Galton-Watson process with offspring distribution $p^{(n)}$ starting from $[nx]$.

\item For $a\in \mbb N$, define the contour function of trees cut at level $a$, $C_a^{p^{(n)},x}=(C_a^{p^{(n)},x}(t), t\geq0)$, by concatenating the contour functions $(C(r_a{\cal G}_{1}^{p^{(n)}}, t), t\in[0,2\#r_a{\cal G}_{1}^{p^{(n)}}]),
     \ldots, (C(r_a{\cal G}_{[nx]}^{p^{(n)}}, t), t\in[0,2\#r_a{\cal G}_{[nx]}^{p^{(n)}}])$ and setting $C_a^{p^{(n)},x}(t)=0$ for $t\geq 2\sum_{i=1}^{[nx]}\#r_a{\cal G}_{i}^{p^{(n)}}$.

\item For $a\in{\mbb R}^+$, define $C_a^{p^{(n)},x}=\pi_a(C_{\lceil a\rceil}^{p^{(n)},x})$.

\item If $\sum_{k\geq0}kp_k^{(n)}\leq 1$, then we define the contour function $C^{p^{(n)},x}=(C^{p^{(n)},x}(t), t\geq0)$  by concatenating the contour functions $(C({\cal G}_{1}^{p^{(n)}}, t), t\in[0,2\#{\cal G}_{1}^{p^{(n)}}]), \ldots, (C({\cal G}_{[nx]}^{p^{(n)}}, t), t\in[0,2\#{\cal G}_{[nx]}^{p^{(n)}}])$ and setting $C^{p^{(n)},x}(t)=0$ for $t\geq 2\sum_{i=1}^{[nx]}\#{\cal G}_{i}^{p^{(n)}}$.
    \end{itemize}

%\section{On convergence of branching processes.}

\noindent
Let $(\gamma_n, n=1,2,\ldots)$ be a nondecreasing sequence of positive numbers converging to $\infty$. Define
$$
G^{(n)}(\lambda)=n\gamma_n[g^{p^{(n)}}(e^{-\lambda/n})-e^{-\lambda/n}],$$  where $g^{p^{(n)}}$ is the generating function of $p^{(n)}$,
and define a probability measure on $[0,\infty)$ by
$$\mu^{(n)}\l(\frac{k-1}{n}\r)=p_k^{(n)},\quad k\geq0.$$
We then present the following statements. By $\overset{(d)}{\rightarrow}$ we mean convergence in distribution.
\begin{enumerate}
\item[(A1)] $G^{(n)}(\lambda)\rar\psi(\lambda)$ as $n\rar\infty$ uniformly on any bounded interval.

\item[(A2)]
\beqlb\label{lem:lib}
\left(\frac{1}{n}Y^{p^{(n)},x}_{[\gamma_n
t]},\; t\geq0\right)\overset{(d)}{\longrightarrow}(Y_t^{\psi,x},\; t\geq0),\quad \text{as }n\rar\infty,
\eeqlb
in $D(\mbb R^+,\mbb R^+)$.

\item[(A3)] There exists a probability measure $\mu$ on $(-\infty, +\infty)$ such that
$\l(\mu^{(n)}\r)^{*[n\gamma_n]}\rar \mu$ as $n\rar\infty$, where $\int e^{-\lambda x}\mu(dx)=e^{\psi(\lambda)}.$
\end{enumerate}
The following lemma is a variant of Theorem 3.4 in \cite{[Gr74]}.
\begin{lemma}\label{lem:li}
Let $\psi$ be a branching mechanism satisfying (H1) and (H2). Then (A1), (A2) and (A3) are equivalent.
\end{lemma}
\begin{remark}
(A3) is just the condition (i) in Theorem 3.4 of \cite{[Gr74]}. Under our assumption on $\psi$, we do not need condition (b) there. (A3) is also equivalent to the convergence of random walks to $\psi$-L\'evy processes; see Theorem 2.1.1 of \cite{[DL02]} for (sub)critical case.
\end{remark}

\proof We shall show that (A2)$\Leftrightarrow$(A3) and (A3)$\Leftrightarrow$(A1).

i): If (A2) holds,  then $\psi$ is conservative implies that ${\mbb P}(Y_t<\infty)=1$ for all $t\geq0$. Then Theorem 3.3 in \cite{[Gr74]} gives (A2)$\Rightarrow$(A3). Meanwhile, Theorem 3.1 in \cite{[Gr74]} implies (A3)$\Rightarrow$(A2).

ii): We first show (A3)$\Rightarrow$(A1). Denote by $L^{(n)}(\lambda)$ the Laplace transform of $\l(\mu^{(n)}\r)^{*[n\gamma_n]}$. Then Theorem 2.1 in \cite{[Gr74]}, together with (A3), gives that for every real number $d>0$
\beqlb\label{lap}
\log L^{(n)}(\lz)=[n\gamma_n] \log \l(\frac{e^{\lz/n}}{n\gamma_n}G^{(n)}(\lz)+1\r)\rar {\psi(\lz)},
\text{ as } n\rar\infty,
\eeqlb uniformly in $\lz\in[0,d],$
which implies that for any $\ez>0$, all $n>n(d,\ez)$ and $\lz\in[0,d]$,
$$
n\gamma_n\l(e^{\frac{\psi(\lz)-\ez}{[n\gamma_n]}}-1\r)
<
e^{\lz/n}G^{(n)}(\lz)<
n\gamma_n\l(e^{\frac{\psi(\lz)+\ez}{[n\gamma_n]}}-1\r).
$$
Then by $|e^x-1-x|<{e^{|x|}}|x|^2/2$,
$$
-\frac{2(\psi(\lz)-\ez)^2 e^{\frac{|\psi(\lz)-\ez|}{[n\gamma_n]}}}{n\gamma_n}-2\ez
<
e^{\lz/n}G^{(n)}(\lz)-\psi(\lz)<
\frac{(\psi(\lz)+\ez)^2 e^{\frac{|\psi(\lz)+\ez|}{n\gamma_n}}}{n\gamma_n}+\ez.$$
Note that $\psi$ is locally bounded. Thus  as $n\rar\infty$,
$G^{(n)}(\lz)\rar \psi(\lz),$
uniformly on any bounded interval, which is just (A1).
Similarly, one can deduce that if $(A1)$ holds, then $ L^{(n)}(\lz)\rar e^{\psi(\lz)}$ as $n\rar\infty$, which implies (A3). \qed

Now, we are ready to present our main theorem. Define $\mathcal{E}^{p^{(n)},x}=\inf\{k\geq0: Y^{p^{(n)},x}_k=0\}$ and $\mathcal{E}^{\psi,x}=\inf\{t\geq0: Y_t^{\psi,x}=0\}$ with the convention that $\inf\emptyset=+\infty$. Denote by $g_{k}^{p^{(n)}}$ the $k$-th iterate of $g^{p^{(n)}}$.
\begin{theorem}\label{Main}
Let $\psi$ be a  branching mechanism satisfying (H1) and (H2).
 Assume that (A1) or (A2) holds.  Suppose in addition that for every $\dz>0$,
\beqlb\label{Main1}
\liminf_{n\rar\infty}g_{[\dz\gamma_n]}^{p^{(n)}}(0)^n>0.
\eeqlb
Then for $x>0$,
\beqlb\label{conexeb}
\frac{1}{\gamma_n}{\cal E}^{p^{(n)},x}\overset{(d)}{\rar}\mathcal{E}^{\psi,x}
\text { on }[0,+\infty]\eeqlb
 and  for any bounded continuous function $F$ on $C({\mbb R}^+, {\mbb R}^+)$ and every $a\geq0$,
\beqlb\label{Main2}
\lim_{n\rar\infty}{\mE}
\l[F\l(\pi_a\l(\gamma_n^{-1}C^{p^{(n)},x}({2n\gamma_n\cdot})\r)\r) \r]=\mE^{\psi,a}_x\l[F\l(e\r)\r].\eeqlb
\end{theorem}

Before proving the theorem, we would like to give some remarks.
\begin{remark}
(\ref{Main1}) is essential to (\ref{Main2}); see the comments following Theorem 2.3.1 in \cite{[DL02]}. In fact under our assumptions $(H1), (H2)$ and $(A1)$, (\ref{Main1}) is equivalent to (\ref{conexeb}). To see (\ref{conexeb}) implies (\ref{Main1}), note that
 $$g_{[\dz\gamma_n]}^{p^{(n)}}(0)^{[nx]}=\mP\l[Y^{p^{(n)},x}_{[\dz \gamma_n]}=0\r]=\mP[{\cal E}^{p^{(n)},x}/{\gamma_n}< \dz]$$
which, together with (\ref{conexeb}), gives
$$
\liminf_{n\rar\infty}g_{[\dz\gamma_n]}^{p^{(n)}}(0)^{[nx]}= \liminf_{n\rar\infty}
\mP[{\cal E}^{p^{(n)},x}/{\gamma_n}<\dz]\geq \mP[{\cal E}^{\psi,x}<\dz]>0,
$$
where the last inequality follows from our assumption $(H1)$; see Chapter 10 in \cite{[Ky06]} for details.
\end{remark}
\begin{remark}\label{remDW}
Some related work on the convergence of discrete Galton-Watson trees has been done in \cite{[DL02]} and \cite{[DW12]}.  In \cite{[DL02]}, only the (sub)critical case was considered; see Theorem \ref{lem:dl} below. In Theorem 4.15 of \cite{[DW12]}, a similar work was done using a quite different formalism. The assumptions there are same as our assumptions in the Theorem \ref{Main}. But the convergence holds for locally compact rooted real trees in the sense of the pointed Gromov-Hausdorff distance, which is a weaker convergence. Thus Theorem \ref{Main} implies that the super-critical L\'evy trees constructed in \cite{ad:ctvmp} coincides with the one studied in  \cite{[DW12]}; see also \cite{[ADH12]} and \cite{[DW07]}.
\end{remark}

We then present a variant of Theorem 2.3.1 and Corollary 2.5.1 in \cite{[DL02]} which is essential to our proof of Theorem \ref{Main}.

\begin{theorem}\label{lem:dl} (Theorem 2.3.1 and Corollary 2.5.1 of \cite{[DL02]})
Let $\psi$ be a (sub)critical branching mechanism satisfying $(H1)$.
 Assume that (A1) or (A2) holds. Suppose in addition that for every $\dz>0$,
\beqlb\label{lem:dlb}
\liminf_{n\rar\infty}g_{[\dz\gamma_n]}^{p^{(n)}}(0)^n>0.
\eeqlb
Then \beqlb\label{lem:dla}
\frac{1}{\gamma_n}{\cal E}^{p^{(n)},x}\overset{(d)}{\rar}\mathcal{E}^{\psi,x}
\text { on }[0,+\infty)\eeqlb
 and for any bounded continuous function $F$ on $C({\mbb R}^+, {\mbb R}^+)\times D({\mbb R}^+, {\mbb R}^+)$,
\beqlb\label{coropia}
\lim_{n\rar\infty}{\mE}
\l[F\l(\pi_a\l(\gamma_n^{-1}C^{p^{(n)},x}({2n\gamma_n\cdot})\r), \l(\frac{1}{n}Y^{p^{(n)}, x}_{[\gamma_n
a]}\r)_{a\geq0}\r) \r]=\mE^{\psi}_x\l[F\l(\pi_a(e), (Z_a)_{a\geq0}\r)\r].
\eeqlb
\end{theorem}
\proof The comments following Theorem 2.3.1 in \cite{[DL02]} give (\ref{lem:dla}). And by Corollary 2.5.1 in \cite{[DL02]}, we have
 \beqlb\label{lem:dlc}
\lim_{n\rar\infty}{\mE}
\l[F\l(\l(\gamma_n^{-1}C^{p^{(n)},x}({2n\gamma_n\cdot})\r), \l(\frac{1}{n}Y^{p^{(n)},x}_{[\gamma_n
a]}\r)_{a\geq0}\r) \r]=\mE^{\psi}_x\l[F\l(e, (Z_a)_{a\geq0}\r)\r].
\eeqlb
 On the other hand, let ${\cal C}_{a}$ be the set of discontinuities of $\pi_a$. (\ref{localtime}) yields
\beqlb\label{coropib}
\Gamma_{e,a}(x)=\int_0^x1_{\{e_t\leq a\}}dt=\int_{\mbb R^+}1_{\{s\leq a\}}L_x^sds=\int_0^x1_{\{e_t<a\}}dt,\quad \mP_x^{\psi}-a.s.
\eeqlb
Then by arguments on page 746 in \cite{[LeG10]},  $\mP_x^{\psi}({\cal D}_a)=0$. Then (\ref{coropia}) follows readily from Theorem 2.7 in \cite{[Bi99]}. \qed

 Recall that $\gamma$ is the largest root of $\psi(\lz)=0$.

\begin{lemma}\label{lem:exe}
Let $\psi$ be a branching mechanism satisfying $(H1)$ and $(H2)$.
 Assume that (A1) or (A2) holds. Suppose in addition that for every $\dz>0$,
\beqlb\label{lem:exe1}
\liminf_{n\rar\infty}g_{[\dz\gamma_n]}^{p^{(n)}}(0)^n>0.
\eeqlb
Then as $n\rar\infty$,
\beqlb\label{conexea}
f(p^{(n)})^{[nx]}\rar e^{-\gamma x},\quad x>0.
\eeqlb
\end{lemma}
\proof
Recall that $f(p^{(n)})$ denotes the minimal solution of $g^{p^{(n)}}(s)=s.$
For each $n\geq1$, define
$$
q^{(n)}_k=f({p^{(n)}})^{k-1}p^{(n)}_k,\quad k\geq1\quad \text{ and }\quad
q^{(n)}_0=1-\sum_{k\geq1}q^{(n)}_k.
$$
Then $q^{(n)}=\{q_k^{(n)}: k\geq0\}$ is a probability distribution with generating function given by
\beqlb\label{gq}
g^{q^{(n)}}(s)=g^{p^{(n)}}\l(sf({p^{(n)}})\r)/f({p^{(n)}}),\quad 0\leq s\leq 1.
\eeqlb

\noindent Thus $ g^{q^{(n)}}(0)=g^{p^{(n)}}(0)/f({p^{(n)}})$ and by induction we further have
\beqlb\label{extinclim}
g_{k+1}^{q^{(n)}}(0)=g^{q^{(n)}}\l(g_{k}^{q^{(n)}}(0)\r)
=g^{p^{(n)}}\l(g_{k}^{q^{(n)}}(0)f({p^{(n)}})\r)/f({p^{(n)}})
=g_{k+1}^{p^{(n)}}(0)/f({p^{(n)}}),\quad k\geq1.
\eeqlb
With (\ref{lem:exe1}),  we see that for any $\delta>0$,
\beqlb\label{sub2}
1\geq\liminf_{n\rar\infty}g_{[\dz \gamma_n]}^{q^{(n)}}(0)^n
=\liminf_{n\rar\infty}g_{[\dz \gamma_n]}^{p^{(n)}}(0)^n/f({p^{(n)}})^n\geq\liminf_{n\rar\infty}g_{[\dz \gamma_n]}^{p^{(n)}}(0)^n>0.
\eeqlb
Then we also have $e^{-\gamma_0}:=\liminf_{n\rar\infty}f({p^{(n)}})^n>0.$
Since $f({p^{(n)}})\leq1$, we may write $f({p^{(n)}})=e^{-a_n/n}$ for some $a_n\geq0$. We further have $\limsup_{n\rar\infty}a_n=\gamma_0.$ We shall show that $\gamma_0=\gamma$ and $\{a_n:n\geq1\}$ is a convergent sequence. To this end, let $\{a_{n_k}:k\geq1\}$ be a convergent subsequence of $\{a_n:n\geq1\}$ with $\lim_{k\rar\infty}a_{n_k}=:\tilde{\gamma}\leq \gamma_0.$ Then by (A1),
$$
0=n_k\gamma_{n_k}[{g^{p^{(n_k)}}(e^{-a_{n_k}/n_k})-e^{-a_{n_k}/{n_k}}}]\rar \psi(\tilde{\gamma}),\quad \text{as}\quad k\rar\infty.
$$
Thus $\psi(\tilde{\gamma})=0.$ On the other hand, note that $\psi$ is a convex function with $\psi(0)=0$ and $\gamma$ is the largest root of $\psi(\lz)=0$. Then we have  $\psi(\lz_1)<0$ and $\psi(\lz_2)>0$ for $0<\lz_1<\gamma<\lz_2$.  If $\tilde{\gamma}\neq\gamma$, then $\tilde{\gamma}=0$. In this case, we may find a sequence $\{b_{n_k}: k\geq1\}$ with $b_{n_k}>a_{n_k}$ for all $k\geq 1$ such that $b_{n_k}\rar\gamma$ and for $k$ sufficiently large
$$
{g^{p^{(n_k)}}(e^{-b_{n_k}/n_k})-e^{-b_{n_k}/{n_k}}}=0.
$$
This contradicts the fact that $f(p^{(n)})=e^{-a_n/n}$ is the minimal solution of $g^{p^{(n)}}(s)=s.$ Thus $\tilde{\gamma}=\gamma$ which implies that $\lim_{n\rar\infty}a_n=\gamma$ and $\lim_{n\rar\infty}f(p^{(n)})^{[nx]}=e^{\gamma x}$ for any $x>0.$  \qed

%%%%%%&&&&&&&&&&&&&&&&&&&&&&&&&&&&&&&&&&&

\bigskip

We are in the position to prove Theorem \ref{Main}.

\bigskip

{\bf Proof of Theorem \ref{Main}:}   With Theorem \ref{lem:dl} in hand, we only need to prove the result when $\psi$ is super-critical. The proof will be divided into three steps.

\textit{First step:} One can deduce from (A1) and (\ref{conexea}) that
\beqlb\label{sub1}
&&n\gamma_n[g^{q^{(n)}}\l(e^{-\lambda/n}\r)-e^{-\lambda/n}]\cr
&&\qquad=n\gamma_n\l
[g^{p^{(n)}}\l(e^{-\lambda/n}f({p^{(n)}})\r)
-e^{-\lambda/n}f({p^{(n)}})\r]/f({p^{(n)}})\cr
&&\qquad\rightarrow\psi(\lambda+\gamma),
\quad
\text{as }n\rightarrow\infty,
\eeqlb
uniformly on any bounded interval. Then Lemma \ref{lem:li} and Theorem \ref{lem:dl}, together with (\ref{sub2}) and (\ref{sub1}), imply that
\beqlb\label{conexesub}
\frac{1}{\gamma_n}{\cal E}^{q^{(n)},x}\overset{(d)}{\rar}\mathcal{E}^{\psi_{\gamma},x}
\text { on }[0,+\infty)\eeqlb
and for any bounded continuous function $F$ on $C({\mbb R}^+, {\mbb R}^+)\times D({\mbb R}^+, {\mbb R})$,
\beqlb\label{consub}
\lim_{n\rar\infty}{\mE}\l[F\l(\pi_a\l((\gamma_n^{-1}C^{q^{(n)},x}({2n\gamma_n\cdot})\r), \l(\frac{1}{n}Y^{q^{(n)},x}_{[\gamma_n
a]}\r)_{a\geq0}\r) \r]=\mE^{\psi_{\gamma}}_x\l[F(\pi_a(e), (Z_a)_{a\geq0})\r].
\eeqlb

\textit{Second step:} We shall prove (\ref{conexeb}). Note that
$$
\{{\cal E}^{p^{(n)},x}<\infty\}=\{{\cal G}_i^{p^{(n)}}, i=1,\ldots,[nx] \text{ are finite trees }\}.
$$
Then by Corollary \ref{cor:con}, for $f\in C({\mbb R}^+,{\mbb R}^+)$,
$$
{\mE}\l[f\l({\cal E}^{p^{(n)},x}/{\gamma_n}\r)1_{\{{\cal E}^{p^{(n)},x}<\infty\}}\r]=f(p^{(n)})^{[nx]}{\mE}\l[f\l({\cal E}^{q^{(n)},x}/{\gamma_n}\r)\r]
$$
which, by (\ref{conexea}), (\ref{conexesub}) and Lemma \ref{AD2.2}, converges to $e^{-\gamma x}{\mE}\l[f\l(\mathcal{E}^{\psi_{\gamma},x}\r)\r]
={\mE}\l[f\l(\mathcal{E}^{\psi,x}\r)1_{\{{\cal E}^{\psi,x}<\infty\}}\r]$, as $n\rar\infty$.
We also have that
$$
\mP[{\cal E}^{p^{(n)},x}=\infty]=1-f(p^{(n)})^{[nx]}\rar1-e^{-\gamma x}=\mP[{\cal E}^{\psi,x}=\infty],\quad\text{as }n\rar\infty,
$$
which gives (\ref{conexeb}).

\textit{Third step:} We shall prove (\ref{Main2}). By Corollary \ref{cor:con}, for any nonnegative measurable function $F$ on $C({\mbb R}^+, {\mbb R}^+)$ and $a\geq0$,
\beqlb\label{defsuper}
{\mE}\l[F(C_{\lc a\rc}^{p^{(n)},x}(\cdot))\r]={\mE}\left[f(p^{(n)})^{[nx]-Y_{\lc a\rc}^{q^{(n)},x}}
F(C_{\lc a\rc}^{q^{(n)},x}(\cdot))\right].
\eeqlb
Note that
$$
C_a^{q^{(n)},x}=\pi_aC_{\lc a\rc}^{q^{(n)},x}\text{ and }\pi_a(\gamma_n^{-1}C^{q^{(n)},x})=\gamma_n^{-1}C_{\gamma_n a}^{q^{(n)},x}.
$$
Then by (\ref{defsuper}) we have for $a\in {\mbb R}^+$
\beqlb\label{defsuper1}
{\mE}\l[F(C_{a}^{p^{(n)},x}(\cdot))\r]={\mE}\left[f(p^{(n)})^{[nx]-Y_{\lc a\rc}^{q^{(n)},x}}
F(C_{a}^{q^{(n)},x}(\cdot))\right]
\eeqlb
and
\beqnn{\mE}\l[F\l(\gamma_n^{-1}
C_{\gamma_na}^{p^{(n)},x}({2n\gamma_n\cdot})\r) \r]
={\mE}\l[f(p^{(n)})^{[nx]-Y_{\lc \gamma_na\rc}^{q^{(n)},x}}F\l(\pi_a\l(\gamma_n^{-1}
C^{q^{(n)},x}({2n\gamma_n\cdot})\r)\r) \r].
\eeqnn
We shall show that $\{f(p^{(n)})^{[nx]-Y_{\lc \gamma_na\rc}^{q^{(n)},x}}, n\geq1\}$ is uniformly integrable.  Write $Y_a^{n}=Y_{\lc \gamma_na\rc}^{q^{(n)},x}/n$ for simplicity.
First, note that ${\mE}\l[f(p^{(n)})^{[nx]-nY_a^{n}}\r]=1$. Then with (\ref{conexea}) and (\ref{consub}) in hand, by the bounded convergence theorem, we have
$$
 \lim_{l\rar\infty}\lim_{n\rar\infty}{\mE}\l[f(p^{(n)})^{[nx]-n(l\wedge Y_a^{n})}\r]
=\lim_{l\rar\infty}{\mE}_x^{\psi_{\gamma}}\l[e^{-\gamma x+\gamma(l\wedge Z_a)}\r]={\mE}_x^{\psi_{\gamma}}\l[e^{-\gamma x+\gamma Z_a}\r]=1.
$$
Note that both ${\mE}_x^{\psi_{\gamma}}\l[e^{-\gamma x+\gamma(l\wedge Z_a)}\r]$ and ${\mE}\l[f(p^{(n)})^{[nx]-n(l\wedge Y_a^{n})}\r]$ are increasing in $l$. Thus for every $\ez>0$, there exist $l_0$ and $n_0$ such that for all $l>l_0$ and $n>n_0$,
$$1-\ez/2<{\mE}\l[f(p^{(n)})^{[nx]-n(l\wedge Y_a^{n})}\r]\leq 1.$$
Meanwhile, since
$$
 \lim_{l\rar\infty}{\mE}\l[f(p^{(n)})^{[nx]-n(l\wedge Y_a^{n})}\r]={\mE}\l[f(p^{(n)})^{[nx]-nY_a^{n}}\r]=1,
$$
there exists $l_1>0$ such that for all $n\geq1$,
$$1-\ez/2<{\mE}\l[f(p^{(n)})^{[nx]-n(l_1\wedge Y_a^{n})}\r]\leq {\mE}\l[f(p^{(n)})^{[nx]-nY_a^{n}}\r]=1.$$
Then for all $n\geq1$,
$${\mE}\l[f(p^{(n)})^{[nx]-nY_a^{n}}1_{\{Y_a^{n}>l_1\}}\r]-
{\mE}\l[f(p^{(n)})^{[nx]-nl_1}1_{\{Y_a^{n}>l_1\}}\r]<\ez/2.$$
Define $C_0=\sup_{n\geq 1}f(p^{(n)})^{[nx]-nl_1}<\infty.$ Then for any set $A\in {\cal F}$ with ${\mP}(A)<\frac{\ez}{2C_0},$
\beqnn
{\mE}\l[f(p^{(n)})^{[nx]-nY_a^{n}}1_A\r]<{\mE}\l[f(p^{(n)})^{[nx]-n(l_1\wedge Y_a^{n})}1_A\r]+\ez/2<\ez.
\eeqnn
Thus $\{f(p^{(n)})^{[nx]-Y_{\lc \gamma_na\rc}^{q^{(n)},x}}, n\geq1\}$ is uniformly integrable; see Lemma 4.10 in \cite{[Ka02]}. Using the Skorohod representation theorem and (\ref{consub}), one can deduce that
\beqnn
&&\lim_{n\rar\infty}{\mE}\l[F\l(\gamma_n^{-1}
C_{\gamma_na}^{p^{(n)},x}({2n\gamma_n\cdot})\r) \r]
\cr&&
\quad=\lim_{n\rar\infty}{\mE}\l[f(p^{(n)})^{[nx]-Y_{\lc \gamma_na\rc}^{q^{(n)},x}}F\l(\pi_a\l(\gamma_n^{-1}
C^{q^{(n)},x}({2n\gamma_n\cdot})\r)\r) \r]
\cr&&
\quad=\mE^{\psi_{\gamma}}_x\l[e^{-\gamma x+\gamma Z_a}F(\pi_a e) \r].
\eeqnn
which is just the right hand side of (\ref{Main2}).
We have completed the proof.
\qed

\begin{remark}Write $C^n_t=\gamma_n^{-1}C^{p^{(n)},x}({2n\gamma_nt})$ for simplicity and  recall that  $(w^a,a\geq0)$ denotes the canonical process on $\cal W$.
Suppose that the assumptions of Theorem \ref{Main} are satisfied. Then one can construct a sequence of probability measures $\bar{\mP}_x^{p^{n}}$  on $\W$ such that for every
$a\geq0$, the distribution of $w^a$ under $\bar{\mP}_x^{p^{n}}$  is the same as $\pi_a(C^n)$ and for $0\leq
a\leq b$,\;
$$
\pi_a(w^b)=\pi_a\quad \bar{\mP}_x^{p^{(n)}}-a.s.
$$
We then have \beqlb\label{conversuper}\bar{\mP}_x^{p^{n}}\rar \bar{\mP}_x^{\psi}\text{ as }n\rar\infty.\eeqlb
\end{remark}

%\section{Convergence towards excursion measure}
\begin{remark}
In \cite{ad:ctvmp}, an excursion measure (`distribution' of a single tree) was also defined. However, we could not find an easy proof of convergence of trees under such excursion measure.
\end{remark}

{\bf Acknowledgement} Both authors would like to give their sincere thanks to J.-F. Delmas and M. Winkel for their valuable comments and suggestions on an earlier version of this paper. H. He is supported by SRFDP (20110003120003), Ministry of Education (985 project) and NSFC (11071021, 11126037, 11201030).
N. Luan is supported by UIBE (11QD17) and NSFC (11201068).

\bigskip\bigskip

\noindent{\small Hui He: Laboratory of Mathematics and
Complex Systems, School of Mathematical Sciences, Beijing Normal
University, Beijing 100875, People's Republic of China. \\
\textit{E-mail:} {hehui@bnu.edu.cn}

\bigskip\bigskip

\noindent{\small Nana Luan: School of Insurance and Economics, University of
International Business and Economics,
Beijing 100029, People's Republic of China. \\
\textit{E-mail:} {luannana@uibe.edu.cn}

\end{document}